\documentclass[12pt]{article}
\usepackage{amsmath}
\usepackage{amssymb}
\usepackage{amscd}
\usepackage{amsthm}

\newtheorem{dfn}{Definition}[section]
\newtheorem{rem}[dfn]{Remark}
\newtheorem{thm}[dfn]{Theorem}

\newtheorem{defn}[dfn]{Definition}

\newtheorem{lem}[dfn]{Lemma}

\newtheorem{prob}[dfn]{Problem}
\newtheorem{prop}[dfn]{Proposition}

\newtheorem{cor}[dfn]{Corollary}
\newtheorem{conj}[dfn]{Conjecture}

\newcommand{\ext}{\operatorname{ext}}
\newcommand{\Exc}{\operatorname{Exc}}
\newcommand{\diam}{\operatorname{diam}}
\newcommand{\Isom}{\operatorname{Isom}}

\def\R{{\mathbb R}}
\def\eps{\epsilon}
\def\al{\alpha}
\def\be{\beta}
\def\ga{\gamma}
\def\Z{\mathbb Z}
\def\V{\mathbb V}
\def\W{\mathbb W}
\def\H{\mathbb H}
\def\N{\mathbb N}
\def\Ga{\Gamma}
\def\del{\delta}
\def\Del{\Delta}
\def\La{\Lambda}

\def\ol{\overline}

\def\si{\sigma}

\def\E{{\mathbb E}}
\def\U{{\mathcal U}}
\def\D{{\mathcal D}}

\def\om{\omega}
\def\Om{{\Omega}}

\def\la{\lambda}
\def\t{\tilde}

\begin{document}

\title{Homological dimension and critical exponent of Kleinian groups}
\author{Michael Kapovich}

\maketitle

\begin{abstract}
We prove an inequality between the relative homological dimension
of a Kleinian group $\Ga\subset \Isom(\H^n)$ and its critical
exponent. As an application of this result we show that for a
geometrically finite Kleinian group $\Ga$, if the topological
dimension of the limit set of $\Ga$ equals its Hausdorff
dimension, then the limit set is a round sphere.
\end{abstract}

\section{Introduction}

One of the frequent themes in the theory of Kleinian groups is
establishing a relation between the abstract algebraic properties
of a Kleinian group and its geometric properties, determined by
its action on the hyperbolic space. Ahlfors finiteness theorem and
Mostow rigidity theorem are among the most important examples of
such relation. In this paper we establish a relation between two
 invariants of a Kleinian group: Virtual homological dimension (an
 algebraic invariant) and the critical exponent (a geometric invariant). We refer the
reader to Section \ref{prelim} for the precise definitions.

Given a Kleinian group $\Ga\subset \Isom(\H^n)$, consider the set
${\mathcal P}$ of its maximal virtually abelian subgroups of
virtual rank $\ge 2$, i.e. the elements of ${\mathcal P}$ are
maximal subgroups which contain a subgroup isomorphic to $\Z^2$.
Form the maximal subset
$$
\Pi:= \{\Pi_i, i\in I\}\subset {\mathcal P}
$$
of pairwise nonconjugate elements of ${\mathcal P}$. In other
words, $\Pi$ consists of representatives of cusps of rank $\ge 2$
in $\Ga$.

We let $vhd_R(\Ga,\Pi)$ and $vcd_R(\Ga,\Pi)$ denote the {\em
virtual homological} and {\em cohomological dimension} of $\Ga$
relative to $\Pi$, where $R$ is a commutative ring with a unit.
(Instead of working with virtual dimensions, one can use the
(co)homological dimension with respect to fields of zero
characteristic, or, more generally, rings where the order of every
finite subgroup of $\Ga$ is invertible.) Let $\del(\Ga)$ be the
{\em critical exponent} of $\Ga$.

Our main result is

\begin{thm}\label{main}
Suppose that $\Ga$ is a virtually torsion-free Kleinian group. Then
$$
vhd_R(\Ga, \Pi)-1\le \del(\Ga).
$$
\end{thm}

\begin{cor}\label{maincor}
Suppose that the pair $(\Ga, \Pi)$ has {\em finite type}, e.g.
$\Ga$ admits a finite $K(\Ga, 1)$ and the set $\Pi$ is finite.
Then
$$
cd_R(\Ga, \Pi)-1\le \del(\Ga).
$$
\end{cor}

One, therefore, can regard these results as either nontrivial
lower bounds on the critical exponent, or as vanishing theorems
for relative (co)homology groups of $\Ga$ with arbitrary twisted
coefficients. These results also can be viewed as generalizing the
classical inequality
$$
\dim(Z)\le \dim_H(Z)
$$
for compact metric spaces $Z$, see \cite{Hurewitz-Wallman}. Here
$\dim(Z)$ is the topological dimension and $\dim_H(Z)$ is the
Hausdorff dimension.

As an application of Corollary \ref{maincor} we prove

\begin{thm}
\label{lattice} Suppose that $\Ga\subset \Isom(\H^n)$ is a
nonelementary geometrically finite group so that the Hausdorff
dimension of its limit set equals its topological dimension $d$.
Then the limit set of $\Ga$ is a round $d$-sphere, i.e. $\Ga$
preserves a $d+1$-dimensional subspace $H\subset \H^n$ and $H/\Ga$
has finite volume.
\end{thm}

This theorem was first proved by Rufus Bowen \cite{Bowen} for
convex-cocompact quasi-fuchsian subgroups of $\Isom(\H^3)$.
Bowen's theorem was extended by Bishop and Jones
\cite{Bishop-Jones} to subgroups of $\Isom(\H^3)$ with parabolic
elements. Bowen's result was generalized by Chenbo Yue \cite{Yue}
to convex-cocompact subgroups of $\Isom(\H^n)$ whose limit sets
are topological spheres, although his argument did not need the
latter assumption. Note that the arguments of Yue do not work in
the presence of parabolic elements. For cocompact discrete groups
of isometries of $CAT(-1)$ spaces, an analogue of Theorem
\ref{lattice} was proved by Bonk and Kleiner
\cite{Bonk-Kleiner01}, see also the work of Besson, Gallot and
Courtois \cite{BCG05}. The latter paper was the inspiration for
our work.

\begin{conj}\label{C3}
Suppose that $\Ga$ is a finitely-generated Kleinian group in
$\Isom(\H^n)$. Then:

1. $d=vcd_R(\Ga, \Pi)-1\le \delta(\Ga)$.

2. In the case of equality, $\Ga$ is geometrically finite and its
the limit set is a round $d$-sphere in $S^{n-1}$.
\end{conj}

Another application of our main theorem is the following property of groups with small critical exponent:

\begin{cor}\label{free}
Suppose that $\del(\Ga)<1$ and $\Ga$ is of type $FP_2$, e.g., is
finitely-presented. Then $\Ga$ is virtually free.
\end{cor}

\begin{prob}
(Cf. Theorem 1.3 in \cite{Bishop-Jones}.) Is it true that every
finitely-generated Kleinian group $\Ga$ with $\del(\Ga)<1$ is
geometrically finite? Is it true that such group is a classical
Schottky-type group?
\end{prob}

The proofs of our results are generalizations of the proofs due to
Besson,  Courtois and Gallot  in \cite{BCG05}. Our main
contribution in comparison to their paper is treatment of
arbitrary coefficient modules, working with relative homology
groups and handling manifolds whose injectivity radius is not
bounded from below. The most nontrivial technical ingredient of
our paper is existence of the {\em natural maps} introduced in
\cite{BCG05} and their properties established in that paper.

\medskip
In the case of finitely-generated Kleinian subgroups $\Ga\subset
\Isom(\H^3)$, our main theorem easily follows from the well-known
facts about $\Ga$. It suffices to consider the case when $\Ga$ is
torsion-free. If $\del(\Ga)=2$, then Theorem \ref{main} states
that
$$
vhd(\Ga,\Pi)\le 3.
$$
The letter inequality immediately follows from the fact that the
hyperbolic manifold $\H^3/\Ga$ is a 3-dimensional
Eilenberg-MacLane space for $\Ga$. Assume therefore that
$\del(\Ga)<2$. Then it follows from the solution of the Tameness
Conjecture \cite{Agol}, \cite{Calegari-Gabai} (which, in turn,
implies Ahlfors' measure zero conjecture) and \cite{Bishop-Jones},
that $\Ga$ is geometrically finite. Therefore either $\Ga$ is a
Schottky-type group or it contains a finitely-generated
quasi-fuchsian subgroup $\Phi\subset \Ga$, whose limit set is a
topological circle. In the later case,
$$
2\ge vhd_R(\Ga, \Pi)\ge vhd_R(\Phi, \Pi\cap \Phi)=2,
$$
while
$$
\del(\Ga)\ge \del(\Phi)\ge 1.
$$
This implies the inequality
$$
1=2-1=vhd_R(\Ga, \Pi)-1\le 1\le \del(\Ga).
$$
If $\Ga$ is a Schottky-type group, then
$$
\Ga\cong F_k * \Pi_1*...*\Pi_m,
$$
where $\Pi_i\in \Pi$ for $i=1,...,m$. Therefore $vhd_R(\Ga,
\Pi)=1$ and Theorem \ref{main} trivially follows.

\bigskip
{\em Sketch of the proof of Theorem \ref{main}.} Let $\eps$ be a positive number which is smaller than the
Margulis constant $\mu_n$ for $\H^n$.
Let $\del:=\del(\Ga)$. We assume that $\Ga$ is torsion-free.
We sketch the proof under the following assumption:

There exists a {\em thick} triangulation of the hyperbolic
manifold $M=\H^n/\Ga$, i.e. a triangulation $T$ and a number
$L<\infty$, so that every $i$-simplex in $T$ not contained in the
$\eps$-thin part $M_{(0,\eps]}$ of $M$ is $L$-bilipschitz
diffeomorphic to the standard Euclidean $i$-simplex. (Existence of
such triangulation was recently proved by Bill Breslin
\cite{Breslin} for $n=3$.)

Suppose that $hd_R(\Ga,\Pi)>\del+1$. Then for some $q> \del+1$,
there exists a flat bundle $\V$ over the manifold $M$, so that
$$
H_q(M, M_{(0,\eps]}; \V)\ne 0.
$$
Pick a chain $\zeta\in C_q(M; \V)$ which projects to a nonzero
class $[\zeta]$ in $H_q(M, M_{(0,\eps]}; \V)$. We then extend
$\zeta$ to the $\eps$-thin part of $M$, to a locally finite
absolute cycle $\hat\zeta$ of finite volume. Besson,  Courtois and
Gallot in \cite{BCG05} proved existence of a {\em natural} map $F:
M\to M$ which is (properly) homotopic to the identity and
satisfies
$$
vol(F_\#(\hat\zeta))\le   \left( \frac{\del +1}{q}\right)^q vol(\hat\zeta).
$$
Since $q>\del+1$, the locally finite cycle
$\hat\zeta_k:=F^k_\#(\hat\zeta)$ satisfies
$$
\lim_{k\to\infty} vol(\hat\zeta_k)=0.
$$
Then we use the deformation lemma of Federer and Fleming to deform
(for large $k$) the cycle $\hat\zeta_k$ to a locally finite cycle
$\hat\xi_k$ which is supported in the $q-1$-skeleton of $T$ away
from $M_{(0,\eps]}$. Therefore $\hat\xi_k$ determines zero
homology class in $H_q(M, M_{(0,\eps]}; \V)$. Since $F^k$ is
properly homotopic to the identity (with uniform control on the
length of the tracks of the homotopy) we conclude that $[\zeta]$
is trivial as well, which is a contradiction.

Since the existence of a thick triangulation is not proven in
general, we use instead a map $\eta$ from $M$ to a simplicial
complex $X$, which is the nerve of an appropriate cover of $M$.
The map $\eta$ is $L_\kappa$-Lipschitz on the $\kappa$-thick part
of $M$ for every $\kappa>0$. This allows us to do the deformation
arguments in $X$ rather than in $T$. This line of arguments is
borrowed from \cite[\S 5.32]{Gromov00}.

{\bf Acknowledgements.}  This work was partially supported by the
NSF  grant DMS 0405180. Most of this paper was written when the
author was visiting the Max Plank Institute for Mathematics in
Bonn. I am grateful to G\'erard Besson and Gilles Courtois for
sharing with me an early version of \cite{BCG05} and to Leonid
Potyagailo for motivating discussions.

\section{Preliminaries}\label{prelim}

\subsection{Geometric preliminaries}\label{geomprelim}

{\bf Basics of Kleinian groups.} We let $\H^n$ denote the
hyperbolic $n$-space, $S^{n-1}$ the ideal boundary of $\H^n$, and
$\Isom(\H^n)$ the isometry group of $\H^n$. A {\em Kleinian group}
is a discrete isometry group of $\H^n$. The limit set of a
Kleinian group $\Ga$ is denoted $\La(\Ga)$. A Kleinian group $\Ga$
is called {\em elementary} if its limit set contains at most 2
points. A Kleinian group is elementary if and only if it is
virtually abelian. We let
$$
Hull(\La(\Ga))\subset \H^n
$$
denote the convex hull of $\La(\Ga)$ in $\H^n$.

Let $\Ga\subset \Isom(\H^n)$ be a Kleinian group, $x\in \H^n$ be a
point and $\eps$ be a positive real number. Let
$$
\Ga_{x,\eps}\subset \Ga$$
 denote the subgroup generated by the
elements $\ga\in \Ga$ such that
$$
d(x, \ga(x))\le \eps.
$$

Then, according to Kazhdan--Margulis lemma, for every $n$ there is
a constant $\mu_n>0$, called the {\em Margulis constant}, such
that $\Ga_{x,\mu_n}$ is elementary, for every  Kleinian subgroup
$\Ga\subset \Isom(\H^n)$ and every point $x\in \H^n$.

\medskip
{\bf Thick-thin decomposition of hyperbolic manifolds.} For a
point $x$ in a Riemannian manifold $M$ (possibly with convex
boundary) define
$$
InRad_M(x)
$$
to be the injectivity radius of $M$ at $x$. Then the function
$InRad_M$ is 1-Lipschitz, i.e., it satisfies
\begin{equation}\label{inrad}
|InRad_M(x)-InRad_M(x')|\le d(x, x').
\end{equation}

Suppose that $M$ is a metrically complete connected hyperbolic
manifold with convex boundary. Let $\t{M}$ denote the universal
cover of $M$. For $0<\eps<\mu_n$ consider the thick-thin
decomposition
$$
M= M_{(0,\eps]}\cup M_{[\eps,\infty)}.
$$
Here thin part $K=M_{(0,\eps]}$ of $M$ is the closure of the set of points $x\in
M$, such that there exists a homotopically nontrivial loop $\ga_x$
based at $x$, whose length is $< \eps$.

Let $K_i$, $i\in J\subset \N$, denote the connected components of
$K$.

\begin{lem}
Each $K_i$ is covered by a contractible submanifold $\t{K}_i$ in
$\H^n$.
\end{lem}
\proof We identify $\pi_1(K_i)$ with an elementary subgroup
$\Pi_i\subset \Ga$. Then $\t{K}_i=\t{K}_i(\eps)$ is the union
$$
\t{K}_i(\eps) = \bigcup_{ \ga\in \Pi_i\setminus \{1\}}
\t{K}_\eps(\ga),
$$
where
$$
\t{K}_\eps(\ga) =\{ z\in \tilde{M}: d(z, \ga(z))\le \eps\}.
$$
Each $\t{K}_\eps(\ga)$ is convex, since the displacement function
of $\ga$ is convex. Of course, the union of convex sets need not
be convex and $K_i$ is, in general, not convex. We first consider
the case when $\Pi_i$ is a cyclic hyperbolic subgroup. Let $A=A_i$
denote the common axis of the nontrivial elements of $\Pi_i$. Then
$A$ is contained in each $\t{K}_\eps(\ga)$. It follows that
$\t{K}_i:=\t{K}_i(\eps)$ is star-like with respect to every point
of $A$. Therefore $\t{K}_i$ is contractible.

If $\Pi_i$ is parabolic, this argument of course does not apply.
Let $\xi=\xi_i$ denote the fixed point of $\Pi_i$. Then $\t{K}_i$
is star-like with respect to $\xi$. Therefore, every map $f:
S^k\to \t{K}_i(\eps)$ can be homotoped to a map $f_\kappa: S^k\to
\t{K}_i(\kappa)$ along the geodesics asymptotic to $\xi$, where
$\kappa$ and
$$
d(\kappa):=\diam(f_\kappa(S^k))$$ can be chosen arbitrarily small.
Then  $f_\kappa(S^k)$ bounds a ball $f_\kappa(B^{k+1})$ within
$d(\kappa)$ from the image of $f_\kappa$. Thus
$$
f_\kappa(B^{k+1}) \subset \t{K}_i(\kappa + 2 d(\kappa)).
$$
By choosing $\kappa$ so that $\kappa + 2 d(\kappa)<\eps$, we
conclude that $\pi_k(\t{K}_i)=0$ for all $k$. \qed

Therefore each $K_i=K(\Pi_i, 1)$ is an Eilenberg-MacLane space for
its fundamental group $\Pi_i$.

\medskip
{\bf Critical exponent of a Kleinian group.} Let $\Ga\subset
\Isom(\H^n)$ be a Kleinian group. Consider the {\em Poincar\'e
series}
$$
f_s=\sum_{\ga\in \Ga} e^{-s d(\ga(o), o)} ,
$$
where $o\in \H^n$ is a base-point and $d$ is the hyperbolic metric
on $\H^n$. Then the {\em critical exponent} of $\Ga$ is
$$
\del(\Ga)=\inf \{s: f_s<\infty\}.
$$

Critical exponent has several alternative descriptions. Define
$$
N(R):= \# \{ x\in \Ga\cdot o : d(x,o)\le R\}.
$$
Then $\del(\Ga)$ is the {\em rate of exponential growth} of
$N(R)$, i.e.
$$
\del(\Ga)=\limsup_{R\to\infty} \frac{\log( N(R))}{R},
$$
see \cite{Nicholls}. Lastly, the critical exponent can be
interpreted in terms of the geometry of the limit set of $\Ga$.

\begin{thm}\label{nich}
(See \cite{Bishop-Jones, Nicholls, Sullivan(1984),Tukia(1984)}.)
For every Kleinian group $\Ga\subset \Isom(\H^n)$, we have:

1. $$ \del(\Ga)= \dim_H(\La_c(\Ga)).
$$
In particular, if $\Ga$ is geometrically finite,
$\La(\Ga)\setminus \La_c(\Ga)$ is at most countable and we obtain
$$
\del_\Ga= \dim_H(\La_c(\Ga)).
$$
2. If $\Ga$ is geometrically finite then either $\La(\Ga)=S^{n-1}$
or $\del(\Ga)<n-1$.
\end{thm}
\noindent Here $\dim_H$ is the {\em Hausdorff dimension} and
$\La_c(\Ga)\subset S^{n-1}$ is the {\em conical limit set} of
$\Ga$.

Thus the critical exponent of a Kleinian group is easy to estimate
from above:
$$
\del(\Ga)\le n-1.
$$
Estimates from below, however, are nontrivial; our main theorem
provides such a lower bound.

\subsection{Algebraic preliminaries}

In this section we collect various definitions and results of
homological algebra. We refer the reader to \cite{Bieri(1976a)},
\cite{Bieri-Eckmann(1978)} and \cite{Brown(1982)} for the detailed
discussion. For the rest of the paper, we let $R$ be a commutative
ring with a unit denoted $1$. We note that although
\cite{Bieri-Eckmann(1978)} and \cite{Brown(1982)} restrict their
discussion to $R=\Z$, the definitions and facts that we will need
directly generalize to the general commutative rings.

\medskip
{\bf Suggestion to the reader.} For most of the paper, the reader
uncomfortable with homological algebra can think of (co)homology
of $\Ga$  with trivial coefficients and of existence of a finite
$K(\Ga,1)$ instead of the finite type condition for $\Ga$. However in the
proofs of Theorem \ref{lattice} and Corollary \ref{free}, we need
(co)homology with twisted coefficients as well as the general
notion of finite type.

\medskip
 A group $\Ga$ is said to be of {\em finite type}, or $FP$
(over $R$), if there exists a resolution by finitely generated
projective $R\Ga$--modules
$$
0\to P_k\to P_{k-1}\to ... \to P_0\to R\to 0.
$$

For instance, if there exists a finite cell complex $K=K(\Ga,1)$,
then $\Ga$ has finite type for every ring $R$. Every group of
finite type is finitely generated, although it does not have to be
finitely-presented, see \cite{Bestvina-Brady}.

More generally, a group $\Ga$ is said to be of type $FP_k$ (over
$R$),  if there exists a partial resolution by finitely generated
{\em projective} $R\Ga$--modules
$$
P_k\to P_{k-1}\to ... \to P_0\to R\to 0.
$$

A group $\Ga$ is said to have {\em cohomological dimension} $k$ if
$k$ is the least integer such that there exists a resolution by
projective $R\Ga$--modules
$$
0\to P_k\to P_{k-1}\to ... \to P_0\to R\to 0.
$$

\begin{lem}\label{brown}
Suppose that $\Ga$ is of type $FP_k$ and $cd(\Ga)\le k$. Then
$\Ga$ is of type $FP$.
\end{lem}
\proof See  discussion following the proof of Proposition 6.1
 in \cite[Chapter VIII]{Brown(1982)}. \qed

\medskip
A group $\Ga$ is said to have {\em homological} (or {\em weak})
{\em dimension} $k$ over $R$, if $k$ is the least integer such
that there exists a resolution by {\em flat} $R\Ga$--modules
$$
0\to F_k\to F_{k-1}\to ... \to F_0\to R\to 0.
$$

Thus the (co)homological dimension of $\Ga$ equals the
(projective) flat dimension of the $\Ga$--module $R\Ga$. The
cohomological and homological dimensions of $\Ga$ are denoted by
$cd_R(\Ga)$ and $hd_R(\Ga)$ respectively. One can restate the
definition of (co)homological dimension in terms of vanishing of
(co)homologies of $\Ga$:

\begin{thm}
(See \cite{Bieri(1976a)}.)
$$
cd_R(\Ga)=\sup \{n: \exists \hbox{~an~} R\Ga\hbox{--module~} V
\hbox{~so that~} H^n(\Ga; V)\ne 0\},
$$
$$
hd_R(\Ga)=\sup \{n: \exists \hbox{~an~} R\Ga\hbox{--module~} V
\hbox{~so that~} H_n(\Ga; V)\ne 0\}.
$$
\end{thm}

\begin{thm}\label{Stallings-Dunwoody}
Let $\Ga$ be a torsion-free group such that $cd_R(\Ga)\le 1$. Then
$\Ga$ is free.
\end{thm}

This theorem was originally proven by Stallings
\cite{Stallings(1968)} for finitely-generated groups and $R=\Z$;
his proof was extended by Swan \cite{Swan(1969)} to arbitrary
groups. Finally, Dunwoody \cite{Dunwoody(1979)} proved this
theorem for arbitrary rings.

\begin{rem}
One can weaken the {\em torsion-free} assumption, by restricting
to groups with torsion of bounded order, see
\cite{Dunwoody(1979)}.
\end{rem}

We will need a generalization of these definitions to the relative
case. In what follows we let $\Ga$ be a group and $\Pi$ be a {\em
nonempty} collection of subgroups
$$\Pi:=\{\Pi_i, i\in I\}.$$

Given an $R\Ga$-module $V$, one defines the relative (co)homology groups
$$
H_*(\Ga, \Pi; V), \quad H^*(\Ga, \Pi; V).
$$

Instead of the algebraic definition of (co)homologies with
coefficients in an $R\Ga$-module $V$, we will be using the {\em
topological} interpretation, following \cite[Section
1.5]{Bieri-Eckmann(1978)}. Let $K:=K(\Ga,1)$ be an
Eilenberg-MacLane space for $\Ga$. Let $C_i:=K(\Pi_i, 1)$, $i\in
I$. We assume that the complexes $C_i$ are embedded in $K$, so
that $C_i\cap C_j=\emptyset$ for $i\ne j$. Set
$$
C:=\bigcup_{i\in I} C_i.
$$

Then we will be computing the (co)homologies of the pair $(\Ga,
\Pi)$ using the relative (co)homologies of $(K, C)$. Namely, let $X$
denote the universal cover of $K$. Let $V_R$ be the module $V$,
regarded as an $R$--module. We obtain the trivial (product) sheaf
$\t{\V}$ over $X$ with fibers $V_R$. We will think of this sheaf
as the sheaf of local (horizontal) sections of the product bundle
$\E:=X\times V_R\to X$. By abusing the notation we will identify
bundles and sheafs of their sections. The group $\Ga$ acts on this
sheaf diagonally:
$$
\ga\cdot (x, v)=(\ga(x), \ga\cdot v), \ga\in \Ga.
$$

The bundle $\E$ (and the sheaf $\t{\V}$) project to the space $K$,
to a bundle $\V\to K$ and its sheaf $\V$ of local horizontal
sections. Then we have natural isomorphisms
$$
H^*(\Ga, \Pi; V)\cong H^*(K, C; \V), \quad H_*(\Ga, \Pi; V)\cong
H_*(K, C; \V).
$$

We will mostly work with the relative homology groups $H_*(K,
C;\V)$, which we will think of as the (relative) singular homology
of $K$ (rel. $C$) with coefficients in $\V$. We refer the reader
to \cite{Johnson-Millson} for the precise definition.

The most important example (for us) of this computation of
relative homologies will be when $\Ga$ is a Kleinian group, the
complex $K$ is the hyperbolic manifold $M=\H^n/\Ga$, and the
subcomplex $C$ is a disjoint union of Margulis tubes and cusps in
$M$. More generally, we will consider the case when $K$ is a
metrically complete connected hyperbolic manifold with convex
boundary.

We now return to the general case of group pairs $(\Ga,\Pi)$.

\begin{defn}
The {\em relative} (co)homological dimension of $\Ga$ (rel. $\Pi$) is defined as
$$
cd_R(\Ga, \Pi)=\sup \{n: \exists \hbox{~an~} R\Ga\hbox{--module~} V \hbox{~so that~} H^n(\Ga, \Pi; V)\ne 0\},
$$
$$
hd_R(\Ga, \Pi)=\sup \{n: \exists \hbox{~an~} R\Ga\hbox{--module~} V \hbox{~so that~} H_n(\Ga, \Pi; V)\ne 0\}.
$$
\end{defn}

In the case of $R=\Z$, we will omit the subscript from the
notation for the (co)homological dimension. Set
$$
R\Ga/\Pi:= \oplus_{i\in I} R\Ga/\Pi_i.
$$
We have the augmentation $\eps: R\Ga/\Pi\to R$, given by $\eps(g
\Pi_i):=1$ for all cosets $g \Pi_i$ and all $i$. Following
\cite[Section 1.1]{Bieri-Eckmann(1978)}, we set
$$
\Del:=\Del_{\Ga/\Pi}:= Ker(\eps).
$$
Then (see \cite[Section 1.1]{Bieri-Eckmann(1978)})
$$
H^k(\Ga, \Pi; V)\cong H^{k-1}(\Ga; Hom(\Del, V)),
$$
$$
H_k(\Ga, \Pi; V)\cong H_{k-1}(\Ga; \Del\otimes V).
$$

The cohomological and homological dimensions of $(\Ga, \Pi)$ can
be interpreted as flat and projective dimensions of
$\Del=\Del_{\Ga/\Pi}$ respectively:
\begin{equation}
\label{E1} hd_R(\Ga, \Pi)-1= flat \dim(\Del), \quad cd_R(\Ga,
\Pi)-1= proj \dim(\Del),
\end{equation}
see \cite[Section 4.1]{Bieri-Eckmann(1978)}. For most of the paper
this interpretation of (co)homological dimension will be
unnecessary; the only exceptions are Lemmata \ref{L1} and \ref{L2}
below:

\begin{lem}\label{L1}
$$
hd_{R}(\Ga, \Pi)\le cd_R (\Ga, \Pi)\le hd_R(\Ga, \Pi)+1.
$$
\end{lem}
\proof The absolute case was proved in \cite{Bieri(1976a)}; the
relative case follows from the same arguments as in Bieri's book
using the equation (\ref{E1}). \qed

\medskip
A pair $(\Ga, \Pi)$ is said to have {\em finite type} (over $R$) if:

1. $\Ga$ and each $\Pi_i$ has type $FP$.

2.   The set $I$ is finite.

\medskip
This condition is stronger than the one considered in
\cite[Section 4.1]{Bieri-Eckmann(1978)}. However it will suffice
for our purposes as we are interested in the case where each
$\Pi_i$ is a finitely generated virtually abelian group. Such
groups $\Pi_i$ necessarily have finite type.

Note that there is a free finitely generated Kleinian group
$\Ga\subset \Isom(\H^4)$, so that $\Ga$ contains infinitely many
$\Ga$-conjugacy classes of maximal parabolic subgroups,
\cite{KP2}. It is unknown if every Kleinian group $\Ga$ of finite
type contains only finitely many conjugacy classes of maximal
parabolic subgroups of rank $\ge 2$.

If $\Ga\subset \Isom(\H^n)$ is a geometrically finite Kleinian
group, then it contains only finitely many conjugacy classes of
maximal parabolic subgroups, see \cite{Bowditch(1993b)}.
 Moreover, $\Ga$ has finite type since its
admits a finite $K(\Ga, 1)$, which is the complement to cusps in
the convex core of $\H^n/\Ga$. Therefore in this case $(\Ga, \Pi)$
has finite type.

\begin{lem}\label{L2}
If $(\Ga, \Pi)$ is of finite type, then

1. $cd_R(\Ga,\Pi)=hd_R(\Ga,\Pi)$.

2. $cd_R(\Ga, \Pi)= \sup \{n:  H^n(\Ga, \Pi; R\Ga)\ne 0\}$.
\end{lem}
\proof This theorem was proved in \cite{Bieri(1976a)} (see also
\cite[Chapter VIII, Proposition 6.7]{Brown(1982)}) in the case
when $\Pi=\emptyset$. The same arguments go through in the
relative case.  \qed 

Suppose that $\Ga$ is virtually torsion free, i.e. it contains a
finite-index subgroup $\Ga'\subset \Ga$ which is torsion-free. Let
$\Pi'$ denote the collection of subgroups of $\Ga'$ obtained by
intersecting $\Ga'$ with the elements of $\Pi$. One defines the
{\em virtual} relative (co)homological dimension of $\Ga$ as
$$
vcd_R(\Ga, \Pi)=cd_R(\Ga', \Pi'),
$$
$$
vhd_R(\Ga, \Pi)=hd_R(\Ga', \Pi').
$$
Recall that every finitely-generated Kleinian group is virtually torsion-free by Selberg's lemma.

\section{Volumes of relative cycles}

Let $X$ be either a simplicial complex or a Riemannian manifold,
possibly with convex boundary. In the case when $X$ is a
simplicial complex, we metrize $X$ by identifying each $i$-simplex
in $X$ with the standard Euclidean $i$-simplex in $\R^{i+1}$. Let
$Y\subset X$ be either a subcomplex or a closed submanifold with
piecewise-smooth boundary. Let $\hat\om_q$ be the $q$-volume form
on $X$ induced by piecewise-Euclidean or Riemannian metric on $X$.
Let $\chi$ be the characteristic function of $X\setminus Y$; we
define the relative $q$-volume form $\om_q$ by
$$
\om_q:= \chi\cdot \hat\om_q.
$$

Let $\W\to X$ be a flat bundle whose fibers are copies of an
$R$-module $V_R$. We define the {\em relative volume} $Vol(\zeta,
Y)$ for piecewise-smooth singular $q$-chains $\zeta$ in $C_q(X,
\W)$ as follows. Consider first the case when $\zeta=w\otimes
\si$, where $\si: \Del^q\to X$ is a singular $q$-simplex and $w$
is a (horizontal) section of $\W$ over the support of $\si$. Then
set
$$
Vol(\zeta, Y)= \int_{\Del^q} \si^*(\om_q).
$$
For a general chain
$$
\zeta=\sum_{i=1}^s w_i\otimes \si_i,
$$
set
$$
Vol(\zeta, Y):= \sum_{i=1}^s Vol(w_i\otimes \si_i).
$$
We set $Vol(\zeta):= Vol(\zeta, \emptyset)$. Clearly, the relative
volume descends to a function on $Z_q(X, Y; \W)$. For a relative
homology class $\xi\in H_q(X, Y; \W)$, we define the relative
volume by
$$
Vol(\xi, Y):= \inf \{Vol(\zeta, Y): \xi=[\zeta]\}.
$$

\medskip Note that our definition of relative volume does not take
the coefficients into account. Suppose that $R$ is a {\em normed}
ring with a norm $|\cdot |$ and $V$ is a {\em normed}
$R\Ga$--module, i.e. it admits a norm $|\cdot |$ such that
$$
|r\ga\cdot v|=|r|\cdot |v|, \quad \forall r\in R, \ga\in \Ga,
$$
where $|r|$ is the norm of $r\in R$. For instance, take $V=\R\Ga$
or $V=\R$, the trivial $\R\Ga$-module. If $R$ is a normed ring,
then the normed modules suffice for calculation of the
cohomological dimension of $\Ga$ over $R$, see \cite[Chapter VIII,
Proposition 2.3]{Brown(1982)}.

Then one can define another volume function, which is sensitive to the coefficients:
$$
vol(w\otimes \si, Y)= \int_{\Del^q} |w|\si^*(\om_q).
$$
However, as the rings discussed in this paper are general (for
instance, we allow finite rings $R$), we cannot use this
definition.

\begin{prob}
Is it true that for every group $\Ga$,
$$
cd_{\R}(\Ga)=\sup \{q: \exists \hbox{~~a Banach~~} \R\Ga\hbox{--module~~} V \hbox{~~so that~~}
H^q(\Ga, V)\ne 0\} \quad ?
$$
Here a {\em Banach} $\R\Ga$--module is a normed $\R\Ga$--module
which is complete as a normed vector space. Note that the answer
is unclear even for groups $\Ga$ of finite type, since $\R\Ga$ is
not a Banach space.
\end{prob}

\section{Coning off singular chains}

Let $M$ be a metrically complete hyperbolic $n$-manifold with convex boundary
and $\V\to M$ be a flat bundle whose fibers are isomorphic to the $R$-module $V_R$. Pick $0<\eps\le \mu_n$. For a
singular chain
$$
\si\in C_q(M; \V),
$$
we define its $\eps$-excision
$$
\Exc_\eps(\si)= \si\cap M_{[\eps,\infty)} \in C_q(M_{[\eps,\infty)}; \V)
$$
by excising the open submanifold $M_{(0,\eps)}\subset M$. The main
goal of this section is to define and examine a converse  to this
procedure.

Let
$$
M_{(0,\eps]}= P\cup Q= P_1\cup ...\cup P_s \cup Q_1\cup ... \cup
Q_l,
$$
where $Q$ is the union of compact components (tubes) $Q_i$ of
$M_{(0,\eps]}$ and $P$ is the union of noncompact components
(cusps)  $P_j$.

\subsection{Extension to the tubes}\label{tubes}

Suppose that $K=Q_i\subset M_{(0,\eps]}$ is a component which
retracts to a closed geodesic $c\subset K$. Given a singular
simplex
$$
\si: \Del^q\to K,
$$
we define the extension $\ext(\si)$ of $\si$ to $\Del^q\times
[0,1]$ as follows. For $x\in \Del^q, t\in [0,1]$, let $x':=\si(x)$
and $x''\in c$ be the point nearest  to $x'$. Choose the point
$$
a=\ext \si (x, t)
$$
on the geodesic segment $\ol{x' x''}$, so that
$$
d(\si(x), a)= t d(x', x'').
$$
We triangulate  $\Del^q\times [0,1]$ so that $\ext(\si)$ is a
singular chain. Finally, extend linearly the operator $\ext$  to
the entire $C_*(K; \V)$.

Suppose that
$$
\zeta\in C_{q+1}(M; \V)
$$
is a chain which projects to a relative cycle in $Z_{q+1}(M,
M_{(0,\eps]};\V)$. For each tube $Q_i$ we consider the extension
$$
\zeta_i':=\ext_{Q_i}(\partial \zeta \cap Q_i).
$$

Since each $Q_i$ retracts to a closed geodesic contained in $Q_i$, we obtain

\begin{lem}
For each $q\ge 1$, the extension
$$
\zeta':= \zeta+ \sum_{i=1}^l \zeta_i'
$$
projects to a relative cycle in $Z_{q+1}(M, P; \V)$.
\end{lem}

\subsection{Extension to the cusps}\label{cusps}

Recall that $P\subset M$ is the union of cusps.  Given a chain
$$
\zeta\in C_{q+1}(M; \V)
$$
which projects to a relative cycle in
$$
Z_{q+1}(M, P; \V),
$$
we will define a locally finite absolute cycle $\hat\zeta$, which is an extension of $\zeta$ to the cusps.

\medskip
Let $\Del^m$ be the standard $m$-simplex $[e_0,...,e_m]$ and $\Del^{m-1}$ be its face $[e_1,...,e_m]$. We parameterize the
{\em punctured simplex}
$$
\Del_\circ^m:= \Del^m \setminus \{e_0\}
$$
as follows. Given a point $z\in \Del_\circ^m$, consider the line
segment $\ol{e_0 x}\subset \Del^m$ containing $z$, where $x\in
\Del^{m-1}$. Then
$$
z= tx+ (1-t)x, 0\le t\le 1.
$$
Therefore we give the point $z$ the coordinates $(x, t)$, $x\in \Del^{m-1}, t\in [0,1]$.

Fix a point $\xi\in \partial \H^n$ and consider a piecewise-smooth singular simplex
$$
\si: \Del^{m-1}\to \H^n.
$$
We define the extension
$$
\ext_\xi(\si): \Del_\circ^m\to \H^n
$$
of $\si$ as follows. For the point $y=\si(x)$ consider the geodesic ray
$$
\rho=\rho_{y,\xi}: [0,\infty) \to \H^n
$$
emanating from $y$ and asymptotic to $\xi$. We parameterize $\rho$ with the unit speed and set
$$
\ext_\xi(\si)(x,t):= \rho( -\log(t)).
$$
Then $\ext_\xi(\si)$ is a piecewise-smooth proper map. Given a singular chain which is a linear combination
$$
\si:=\sum_i w_i\otimes \si_i, \quad w_i\in V_R,
$$
we set
$$
\ext_\xi(\si):= \sum_i w_i\otimes \ext_\xi(\si_i).
$$
This extension satisfies
$$
\partial \ext_\xi(\si) = \ext_\xi(\partial \si).
$$
The extension is invariant under the action of $\Isom(\H^n)$ in the sense that
$$
\ga_*(\ext_\xi(\si))= \ext_{\ga(\xi)}(\ga_*(\si)), \quad \forall \ga\in \Isom(\H^n).
$$

For a chain
$$
\si=\sum_i w_i\otimes \si_i,
$$
we define volume of the ``punctured'' chain $\ext_\xi(\si)$ by
$$
Vol(\ext_\xi(\si)):= \sum_i Vol(\ext_\xi(\si_i)).
$$

\begin{lem}
\label{bound}
For every $q$-chain $\si$, $q\ge 1$, we have
$$
Vol(\ext_\xi(\si))\le q\cdot Vol(\si).
$$
\end{lem}
\proof It suffices to prove the inequality in the case of a singular simplex $\si$. We will work in the upper half-space model
of $\H^n$, so that $\xi=\infty$. By subdividing the chain $\si$ appropriately we can assume that
$\si(x)=(x, f(x))$, is the graph of a continuous map
$$
f: \Om \to (0,\infty),
$$
where $\Om$ is a bounded domain in $\R^q\subset \R^{n-1}\subset \partial \H^n$, and $f$ is smooth on the interior of $\Om$.
Then
$$
Vol(\si)= \int_\Om \frac{\sqrt{1+ |\nabla f|^2}}{f(x)^{q}}dx \ge \int_\Om \frac{dx}{f(x)^{q}}
$$
and
$$
Vol(\ext_\xi(\si))= \int_{\Om} \int_{f(x)}^\infty \frac{dt}{t^{q+1}} dx = \int_{\Om}  \frac{q dx}{f(x)^{q}}.
$$
Therefore
$$
Vol(\ext_\xi(\si))\le q\cdot Vol(\si). \qed
$$

Suppose now that $\zeta\in C_{q+1}(M; \V)$,
$$
\partial \zeta=\sum_{j=1}^s \zeta_j, \quad \zeta_j\in C_q(P_j; \V), i=1,...,s.
$$
For every singular chain
$$
\si=\zeta_j=\sum_i w_i \otimes \si_i \in C_q(P_j, \V),
$$
we define a locally finite singular chain
$$
\ext(\si)\in C^{lf}_{q+1}(P_j, \V)
$$
as follows.  Lift each $\si_i$ to a chain
$$
\t\si=\sum_i v_i \otimes \t\si_i \in C_q(\t{P}_j, V_R), v_i\in
V_R,
$$
where $\t{P}_j\subset \H^n$ is a component of the preimage of $P_j$; let $\Pi_j$ be the stabilizer of $\t{P}_j$ in $\Ga$.
Let $\xi=\xi_j$ be a point fixed by $\Pi_j$.

\begin{rem}
Our construction does not depend on whether $\Pi_j$ is parabolic or hyperbolic. In parabolic case we, of course,
have unique fixed point.

If $\Pi_j$ is hyperbolic, then the extension chain $\ext(\si)$ below is not going to be
locally finite in $M$, as it ``spins towards'' a closed geodesic.
\end{rem}

Now extend $\t\si$ to a punctured chain $ext_\xi(\t\si)$. Finally,
project the latter to a punctured chain $\ext(\si)$ via the
universal cover $\H^n\to M$. Triangulate $\Del^{q+1}_\circ$, so
that $\ext(\si)$ is a locally-finite singular chain. Since
$\t{P}_j$ is star-like with respect to $\xi$, it follows that
$$
\ext(\si)\in C^{lf}_{q+1}(P_j, \V).
$$
Invariance of the extension $\ext_\xi$ under $\Isom(\H^n)$, ensures that $\ext(\si)$ does not depend on the choice
of the lifts $\t\si_i$. We also have
$$
\partial \ext(\si) = \ext(\partial \si).
$$

Lemma \ref{bound} implies

\begin{cor}\label{cor:bound}
(Cf. \cite{Thurston(1986)}.) $Vol(\ext(\si))\le q\cdot Vol(\si)$.
\end{cor}

Lastly, set
$$
\hat\zeta:= \zeta+ \sum_{j=1}^s \ext(\zeta_j)
$$

The excision operation obviously extends to the locally-finite chains $\ext(\si)$ and we obtain
$$
\Exc_{\kappa}(\hat\zeta) \in Z_{q+1}(M, M_{(0,\kappa])}; \V)
$$
for every $0<\kappa\le \eps$.
It is clear that
$$
[\Exc_{\kappa}(\hat\zeta)]=[\hat\zeta]\in H_q(M, M_{(0,\eps]}; \V).
$$

Therefore we obtain

\begin{prop}
Let $0<\eps\le \mu_n$, $P\subset M_{(0,\eps]}$ be the union of
cusps. Then for every chain $\zeta \in C_{q+1}(M; \V)$, which
projects to a relative homology class in $H_q(M, P; \V)$, there
exists a locally finite cycle
$$
\hat\zeta\in Z^{lf}_{q+1}(M; \V),
$$
so that:

1. $Vol(\hat\zeta)\le q Vol(\partial \zeta) +Vol(\zeta)$.

2. $[\Exc_{\kappa}(\hat\zeta)]=[\hat\zeta]\in H_q(M, M_{(0,\eps]}; \V)$ for every $0<\kappa\le \eps$.
\end{prop}

\section{Cuspidal homology}

Let $M$ be a metrically complete hyperbolic $n$-manifold with
convex boundary. Let $\pi_1(M)\cong \Ga\subset \Isom(\H^n)$ and
$\Pi$ be the collection of  cusps in $\Ga$, i.e. $\Pi$ consists of
representatives of $\Ga$-conjugacy classes of maximal parabolic
subgroups of $\Ga$. Note that here we allow parabolic subgroups of
rank $1$. Let $\V$ be a flat bundle over $M$ associated with an
$R\Ga$-module $V$. The elements $\Pi_i$ of $\Pi$ correspond to the
components $P_i$ of $P\subset M_{(0,\mu_n]}$. Given $0<\eps\le
\mu_n$, consider the thick-thin decomposition
$$
M=M_{(0,\eps]}\cup M_{[\eps,\infty)}
$$
and let
$$
P_\eps:=P\cap M_{(0,\eps]}.
$$

We then have the direct system
$$
(M, P_\eps)\to (M, P_{\eps'}), \quad 0<\eps\le \eps'\le \mu_n.
$$
These maps induce isomorphisms
$$
H_*(M, P_\eps; \V)\to H_*(M, P_{\eps'}; \V).
$$
We therefore identify
$$
\lim_{\eps} H_*(M, P_\eps; \V)\cong H_*(M, P_{\mu_n}; \V) \cong H^*(\Ga, \Pi; V).
$$
We will refer to this direct limit as the {\em cuspidal homology} of $M$,
$$
H^{cusp}_*(M; \V).
$$
We have an obvious homomorphism
$$
Exc: H^{lf}_*(M; \V)\to H^{cusp}_*(M; \V)
$$
given by the excision. The advantage of working with the above direct limit is the following:

Suppose that $f: M\to M$ is a proper $L$-Lipschitz map. Then
$f$ induces a self-map
$$
f: \{(M, P_\eps)\}\to \{(M, P_\eps)\}
$$
of the direct system; the latter clearly induces an isomorphism
$$
f_*: H^{cusp}_*(M; \V)\to H^{cusp}_*(M; \V).
$$
If $f$ induces the identity automorphism of $\Ga$, then $f_*=Id$.

\section{Partition of unity and a map to the nerve for a hyperbolic $n$-manifold}

Fix a number $\eps>0$. Suppose that $M$ is a complete hyperbolic
$n$-manifold. (In this section we do not allow $M$ to have
boundary.) Given a covering $\U$ of $M$ by contractible open sets
with contractible intersections, $M$ is homotopy--equivalent to
the nerve of $\U$. The goal of this section is to get a
homotopy--equivalence with controlled Lipschitz constant. The
Lipschitz constant will be bounded on the thick part of $M$. This
construction is standard (cf. \cite[\S 5.32]{Gromov00}), we
include it for the sake of completeness.

Our first goal is to find an appropriate covering $\U$ by convex
metric balls $B_{\eps_j}(x_j)$ in $M$, whose radii $\eps_j$ are
multiples of the injectivity radii of $x_j$.

Choose $0<\al< 1$ and define the  function
$$
\tau(x):= \al \min (InRad_M(x), \frac{\eps}{2}).
$$
Note that $\tau$ is a continuous function which is constant on $M_{[\eps,\infty)}$. Recall
that $B_r(x)$ denotes the open $r$-ball centered at $x$.

\begin{lem}\label{compare}
Suppose that $x, y\in M$ are such that
$$
B_{\tau(x)}(x)\cap B_{\tau(y)}(y) \ne \emptyset.
$$
Then
$$
\frac{\tau(x)}{\tau(y)}\le \frac{1+\al}{1-\al}.
$$
\end{lem}
\proof If $x, y\in M_{[\eps,\infty)}$ then $\tau(x)=\tau(y)=\al\eps$ and we are done.
We consider the case
$$
\tau(x)\ge  \tau(y)=\al\cdot InRad_M(y).
$$
Then the inequality (\ref{inrad}) implies that
$$
\al^{-1}(\tau(x)-\tau(y))\le InRad_M(x)- InRad_M(y) \le d(x,y)\le \tau(x)+\tau(y).
$$
Therefore
$$
\frac{\tau(x)}{\tau(y)}\le \frac{1+\al}{1-\al}.     \qed
$$

\begin{lem}
\label{cover}
Let $\be>0$ be such that
$$
0< \frac{\be}{\frac{1}{2}-\be} < \frac{1-\al}{1+\al}.
$$
Then there exists a covering $\D$ of the manifold $M$ by the open balls
$$
D_i:=B_{\tau(x_i)/2}(x_i), i\in I,
$$
so that
\begin{equation}\label{B2}
B_{\be\tau(x_j)}(x_j)\cap B_{\be \tau(x_i)}(x_i)=\emptyset, \quad
\forall ~x_i\ne x_j.
\end{equation}
\end{lem}
\proof We construct the set $E$ of centers $x_i$, $i\in I$, of the above balls as follows. Choose a maximal set
$$
E=\{x_i, i\in I\}\subset M,
$$
satisfying (\ref{B2}). Suppose that $\D:=\{D_j, j\in I\}$ is not a covering of $M$. Then there exists
$$
x\in M \hbox{~~so that~~} d(x, x_i)\ge \tau(x_i)/2, \quad \forall i\in I.
$$
Suppose that there exists $i\in I$ so that
$$
B_{\be\tau(x)}(x)\cap B_{\be \tau(x_i)}(x_i)\ne \emptyset.
$$
Then
$$
\frac{\tau(y)}{2}\le d(x, y)<\be(\tau(x)+\tau(y))
$$
for $y:=x_i$. By Lemma \ref{compare},
$$
\frac{\tau(x)}{\tau(y)}\le \frac{1+\al}{1-\al}
$$
and, by combining these inequalities, we get
$$
\frac{\be}{\frac{1}{2}-\be} \ge \frac{1-\al}{1+\al}.
$$
This contradicts out choice of $\be$. Therefore
$$
B_{\be\tau(x)}(x)\cap B_{\be \tau(x_i)}=\emptyset, \quad \forall i\in I,
$$
which contradicts maximality of $E$. Hence
$$
\D=\{ B_{\tau(x_i)/2}(x_i), x_i\in E\}
$$
is the required covering of $M$. \qed

\medskip
We leave it to the reader to verify that if $0< \al\le 1/8$, $0<\be$, and
$$
0< \frac{\be}{\frac{1}{2}-\be} < \frac{1-\al}{1+\al},
$$
then
\begin{equation}
\label{albe}
1+\be \le \frac{1- \al}{\al(1+\al)}.
\end{equation}

We let $vol(r)$ denote the
volume of a hyperbolic $r$-ball in $\H^n$.
Let $\om_n$ denote the volume of the Euclidean $n$-ball of the unit radius.
Then
\begin{equation}\label{I1}
vol(r)\ge \om_n r^n.
\end{equation}
Moreover, whenever $r\le 1$, we have
\begin{equation}
\label{I2}
vol(r)\le \om_n\cdot 2^{n-1} r^n.
\end{equation}

Given a choice of $\al$ and $\be$ as above, define the covering
$$
\U:= \{ B_i=B_{\tau(x_i)}(x_i)  : i\in I\}
$$
of the manifold $M$, where $E=\{x_i, i\in I\}$ is as in Lemma
\ref{cover}.

\begin{lem}\label{multiplicity}
Assume that $0<\al\le 1/8$,\ $0<\be$ and
$$
0< \frac{\be}{\frac{1}{2}-\be} < \frac{1-\al}{1+\al}.
$$
Then the covering $\U$ has multiplicity at most
$$
\frac{2^{2n-1}}{\be^n}.
$$
\end{lem}
\proof Suppose that
$$
y\in \bigcap_{j\in J} B_j,
$$
for some $J\subset I$ of cardinality $m$. Let $\nu:=\max\{\tau(x_j), j\in J\}$.
Then for every $j\in J$,
$$
B_{\be\tau(x_j)}(x_j)\subset B_{(1+\be)\nu}(y).
$$
The balls $B_{\be\tau(x_j)}(x_j)$, $j\in J$, are pairwise disjoint by the definition of $\U$. Therefore
$$
\sum_{j\in J} Vol(B_{\be\tau(x_j)}(x_j)) \le Vol (B_{(1+\be)\nu}(y)).
$$
Recall that
$$
\frac{\tau(y)}{\al}\le InRad_M(y).
$$
By combining this with the inequality
$$
\nu=\tau(x_j)\le \frac{1+\al}{1-\al}\tau(y)
$$
(for some $j\in J$) and the inequality (\ref{albe}), we obtain
$$
(1+\be) \nu \le InRad_M(y).
$$
Therefore the ball $B_{(1+\be)\nu}(y)$ is contained in a normal
ball and thus has the volume
$$
vol( (1+\be)\nu ).
$$

By Lemma \ref{compare},
$$
\tau(x_i) \ge \nu \frac{1-\al}{1+\al}, \quad \forall i\in J.
$$
Hence
$$
 m \cdot vol\left(\be
\nu \frac{1-\al}{1+\al} \right)\le \sum_{i\in J}
Vol(B_{\be\tau(x_i)}(x_i))\le   Vol (B_{(1+\be)\nu}(y))= vol(
(1+\be)\nu ).
$$
Combining this with the inequalities (\ref{I1}) and (\ref{I2}), we
obtain
$$
m\le 2^{n-1} \left[\frac{ (1+\be)(1+\al)}{\be(1-\al)} \right]^n.
$$
Since $\al\le 1/8$ and $\be<1/2$, it follows that
$$
m \le \frac{2^{2n-1}}{\be^n}. \qed
$$

We now fix $\al>0$ and $\be>0$ such that
$$
\al\le \frac{1}{8}
$$
and
$$
0< \frac{\be}{\frac{1}{2}-\be} < \frac{1-\al}{1+\al}.
$$
For instance, take $\al=\be=1/8$. We then obtain

\begin{prop}\label{covering}
There exists a function $m(n, \eps): \N\times (0,\infty)\to \N$,
with the following property. For every complete hyperbolic
$n$-manifold $M$, there exits a countable subset $E=\{x_i, i\in
I\}\subset M$ and a collection of positive numbers $\{\rho_i, i\in
I\}$, so that:

\smallskip
1. Set $\D:= \{D_i=B_{\rho_i/2}(x_i): x_i\in E\}$ and $\U:=  \{B_i=B_{\rho_i}(x_i): x_i\in E\}$.
Then $\D$ (and therefore $\U$) covers $M$.

2. For every  $x_i\in M_{[\eps,\infty)}$,
$$
\rho_i= \frac{\eps}{16}.
$$

3. For every  $x_i\in M_{(0, \eps)}$,
\begin{equation}\label{inradius}
\frac{InRad_M(x_i)}{8} = \rho_i \le \frac{\eps}{16}.
\end{equation}

4. The multiplicity of the covering $\U$ is at most $m(n, \eps)$.
\end{prop}
\proof Set $\rho_i:= \tau(x_i)$. The rest follows from Lemmata \ref{cover} and  \ref{multiplicity}.
\qed

\begin{cor}\label{normal}
Suppose that
$$
x\in \bigcap_{j\in J} B_j,
$$
where $J\subset I$. Then Part 3 of Proposition \ref{covering} in
conjunction with the inequality (\ref{inrad}), implies that the
union
$$
\bigcup_{j\in J} B_j
$$
is contained in a normal neighborhood $N_x$ of $x$.
\end{cor}

We now associate a partition of unity $\{\eta_i, i\in I\}$ to the
covering $\U$ as in Proposition \ref{covering}. For the metric
$r$-ball $B_r(o)\subset \H^n$ we define a {\em bump-function}
$b_r(x)$ on $\H^n$ supported in $B_r(o)$, so that:

1. $0\le b_r(x)\le 1$, $\forall x\in \H^n$.

2. $b_r(x)=1$ for all $x\in B_{r/2}(o)$.

3. $\|\nabla b_r(x)\| \le \psi(r^{-1})$, $\forall x\in \H^n$, where $\psi(r)$ is a
continuous function on $[0,\infty)$, which vanishes at $0$.

\begin{lem}
\label{partition} There exists a smooth partition of unity
$\{\eta_i, i\in \N\}$ subordinate to the covering $\U$, so that
every function $\eta_i$ is $l_i$-Lipschitz, with
$$
L_\kappa:=\sup \{l_i: x_i\in M_{[\kappa,\infty)} \}<\infty
$$
for every $\kappa>0$.
\end{lem}
\proof For every $i\in I$ consider the bump-function $b_i$ on the ball
$B_i$, which equals to $b_{\rho_i}$ after the isometric identification of $B_i$ with
the ball $B_{\rho_i}(o)\subset \H^n$. Note that the radii of the balls $B_i$ are
at least
$$
\min( \frac{\kappa}{16}, \frac{\eps}{16})
$$
for all
$$
x_i\in M_{[\kappa,\infty)}.
$$
Therefore
$$
\|\nabla b_i(x)\|\le \lambda(\kappa)=\max(\psi(   \frac{16}{\kappa}  ), \psi(\frac{16}{\eps} ) )
$$
for such $x_i$. Define the smooth partition of unity $\{\eta_i, i\in \N\}$ subordinate to
the covering $\U$, using the bump-functions $b_i$:
Set
$$
c(x):= \sum_{j\in I}b_j(x),
$$
and
$$
\eta_i(x):= \frac{b_i(x)}{c(x)}.
$$
Then Parts 1 and 4 of Proposition \ref{covering} imply that for every $x\in M$,
$$
1\le c(x) \le m=m(n, \eps).$$
It follows that for every
$$
x_i\in M_{[\kappa,\infty)}
$$
we have
$$
\|\nabla \eta_{i}(x) \|\le \sqrt{n}(m+1) \lambda(\kappa)=L_\kappa. \qed
$$

\medskip
Let $\{\eta_i, i\in I\}$ be a partition of unity for $M$ as above, subordinate to the covering $\U$ of
multiplicity $m\le m(n, \rho_0)$. We identify $I$ with a subset of $\N$.
Let $X$ denote the $m$-dimensional simplicial complex which is the nerve of the covering $\U$.
Collections of balls $\{B_j, j\in J\}$ such that
$$
\bigcap_{j\in J} B_{j} \ne \emptyset
$$
correspond to simplices $\Del_J\subset X$.

We now use the above partition of unity to define a map $\eta$ from $M$ to $X$.
Set
$$
\Del^\infty:= \{ z\in \R^\infty: \sum_{i=1}^\infty z_i=1, z_i\ge
0, i\in \N\},
$$
where
$$
\R^\infty=\bigoplus_{i\in \N} \R.
$$
Then $X$ embeds naturally in $\Del^\infty$. We define the map
$$
\eta: M\to \Del^\infty
$$
by
$$
\eta(x)=(\eta_1(x),...,\eta_k(x),...).
$$
Since$\{(U_i,\eta_i)\}$ is a partition of unity, it is clear  that the map $\eta$ is well-defined.
Moreover, the image of $\eta$ is contained in $X\subset \Del^\infty$. Lemma  \ref{partition} implies

\begin{cor}\label{Lip}
The map $\eta: M\to X$ is piecewise-smooth and
$$
\eta| M_{[\kappa,\infty)}
$$
is $\sqrt{m}L_\kappa$-Lipschitz for every $\kappa>0$.
\end{cor}

Since $\U$ is a covering by convex sets, the map $\eta: M\to X$ is a homotopy-equivalence. Our goal is to construct
its homotopy-inverse $\bar\eta$ with uniform control on the length of the  tracks of the homotopy.

Let $\Del_J=[e_{j_0},...,e_{j_k}]$ be a $k$-simplex in $X$, $J=\{j_0,...,j_k\}\subset I$.
Recall that the vertices $e_j, j\in J$ correspond to the balls
$B_j=B_{\rho_j}(x_j)\in \U$. We define the map
$$
\bar\eta: \Del_J\to M
$$
by sending the vertices $e_j$ to the corresponding centers $x_j\in B_j$.
The union of the balls
$$
\bigcup_{j\in J} B_j
$$
is contained in a normal neighborhood $N_x$ in $M$ (see Corollary \ref{normal}). Consider the convex hull
$$
Hull_J=Hull(\{x_j, j\in J\})\subset N_x.
$$
Since $N_x$ is a normal neighborhood, we can regard $N_x$ as a subset of $\H^n$.
We use the projective model for $Hull_J$. Then there exists a canonical projective map
$$
\bar\eta: \Del_J\to Hull_J
$$
which extends the map defined on the vertices of $\Del$. Namely, the projectivization
$P: \R^{n+1}\setminus \{0\}\to \R P^n$, identifies
$\H^n$ with the hyperboloid
$$
H^n:=\{(t_0,...,t_n)\in \R^{n+1}: t_0>0, -t_0^2+ t_1^2+...+t_n=-1\}.
$$
The points $x_j$ are projections of the points $\hat{x}_j\in H^n$.
Now, there exists a unique linear map
$$
\Del_J\to \R^{n+1}
$$
which sends each $e_j$ to $\hat{x}_j, j\in J$.
Let $\bar\eta$ be the composition of this map with the projection $P$.

\begin{rem}
This is the only place in our argument where we used the fact that $M$ has constant curvature. One can avoid using
the canonical projective map by appealing to convexity of the balls $B_j$ and defining the map $\bar\eta$ (noncanonically)
by the induction on skeleta of $X$.
\end{rem}

We now estimate the displacement for the composition $\bar\eta\circ \eta$. Set $\kappa:=\eps/8$.

\begin{lem}
1. If $x\in M$ and $z\in Star(\eta(x))$, then $d(x, \bar\eta(z))\le \kappa$.

2. There exists a homotopy $H$ between $\bar\eta\circ \eta$ and $Id$, whose tracks
have length $\le \kappa$.
\end{lem}
\proof 1. Let $\Del_J\subset X$ be the smallest simplex containing $\eta(x)$.
If $z\in Star(\eta(x))$, then $z\in \Del_{J'}$, where $J'\subset I$ is a subset containing $J$, so that
$$
\bigcap_{i\in J'} B_i\ne \emptyset.
$$
Let $y:= \bar\eta(z)$. Then
$$
y\in \bigcup_{j\in J'} B_j.
$$
It follows that
$$
d(x, y)\le 2(\rho_i+\rho_j)
$$
for some $j\in J, i\in J'$. Since $\rho_k\le \eps/16$ for all $k\in I$, the first assertion of lemma follows.

2. Part 1 clearly implies that
$$
d(\bar\eta\circ \eta, Id)\le \kappa=\eps/8.
$$
Moreover, for $x\in M$,
$$
y=\bar\eta \circ \eta(x)\in \bigcup_{j\in J} B_j,
$$
where $J\subset I$ is defined as above. Therefore the points $x$ and $y$ belong to the convex hull
$Hull_J\subset N_x$. Therefore we can take the homotopy $H$ between $\bar\eta\circ \eta$
and $Id$ to be the geodesic homotopy along geodesics contained in the convex sets $Hull_J$,
$J\subset \N$ are such that $\Del_J\subset X$.
The length of the tracks of this homotopy clearly does not exceed $\kappa$. \qed

\bigskip
Define the subcomplex $Y=Y_\eps\subset X$ to be the star of
$\eta(M_{(0,\eps]})$ in $X$. Hence we have the map of pairs
$$
\eta: (M, M_{(0,\eps]}) \to (X, Y_\eps).
$$
Note that, because of the lack of convexity of the components of $M_{(0, \eps]}$,
it is unclear if $\eta$ is a homotopy-equivalence of pairs. Nevertheless,
$\eta: M\to X$ is a homotopy-equivalence, hence it induces an isomorphism of the fundamental groups
$$
\Ga=\pi_1(M) \to \pi_1(X).
$$
We obtain a flat bundle $\W$ over $X$, associated with the $R\Ga$-module $V$.

We will see that the map
$$
\eta_{\eps,*}: H_*(M,  M_{(0,\eps]}, \V) \to H_*(X, Y_\eps; \W)
$$
induced by $\eta$, is an {\em approximate} monomorphism, in the following sense:

\begin{prop}
\label{promono}
The kernel of $\eta_{\eps,*}$ is contained in the kernel of
$$
 H_*(M,  M_{(0,\eps]}, \V) \to  H_*(M,  M_{(0, 2\kappa +\eps]}, \V).
$$
\end{prop}
\proof Let $\al\in C_q(M, \W)$ be such that
$$
[\al]\in Ker(\eta_{\eps,q})\subset H_q(M,  M_{(0,\eps]}, \V).
$$
Let $\be\in C_{q+1}(X, Y; \W)$ be a chain so that
$$
\eta_{\#}(\al) -\partial \beta \in C_q(Y; \W).
$$
Set
$$
\al':= (\bar\eta\circ \eta)_{\#}(\al).
$$
Then
$$
\al'- \partial \bar\eta_{\#}(\beta) \in C_q(\bar\eta(Y); \V)
$$
Since $\bar\eta(Y)$ is contained in the $\kappa$-neighborhood of
$M_{(0,\eps]}$, it follows that
$$
\bar\eta(Y)\subset M_{(0,\eps+ 2\kappa]}.
$$
Therefore
$$
\al'\in B_q(M,  M_{(0,\eps+ 2\kappa]}; \V).
$$
On the other hand, since the tracks of the homotopy $H$ have length $\le \kappa$, it gives us a chain
$$
\be'\in C_{q+1}(M,  M_{(0,\eps+ 2\kappa]}; \V),
$$
so that
$$
\al- \al' -\partial \be' \in C_*( M_{(0,\eps+ 2\kappa]}; \V).
$$
Thus
$$
\al \in B_q(M,  M_{(0,\eps+ 2\kappa]}; \V). \qed
$$

\begin{rem}
One can also prove that $\eta_{\eps,*}$ is an {\em approximate}
epimorphism (see \cite{Kapovich-Kleiner05} for the definition),
but we will not need this.
\end{rem}

\section{Vanishing of relative homology classes of small volume}

The main goal of this section is to prove the following (cf.
\cite[Theorem 5.38]{Gromov00}):

\begin{thm}\label{vanishing}
There exists a function $\theta=\theta_n(\eps)$ with the following property. Let $0<\eps<\mu_n/4$.
Let $M$ be a complete (connected) hyperbolic $n$-manifold
with the fundamental group $\Ga$ and the thick-thin decomposition
$$
M=M_{(0,\eps]}\cup M_{[\eps,\infty)},
$$
where $\mu=\mu_n$ is the Margulis constant. Let $\V\to M$ be the flat bundle associated with a $R\Ga$-module $V$.
Then  every relative homology class
$$
[\zeta]\in H_q(M, M_{(0,\eps]}; \V), q>0,
$$
whose (relative) volume is less than $\theta_n(\eps)$, is trivial.
\end{thm}
\proof Let $\eta: M\to X$ be the homotopy-equivalence from $M$ to the nerve of an appropriate cover,
constructed in the previous section. Let $m$ denote the dimension of $X$. As $\eta$ defines an isomorphism
$$
\pi_1(M)\to \pi_1(X),
$$
we obtain the flat bundle $\W$ over $X$, associated with the $R\Ga$-module $V$.
The map $\eta$ induces an approximate isomorphism
$$
\eta_*: H_q(M, M_{(0,\eps]}; \V) \to H_q(X, Y_\eps; \W),
$$
where $Y:=Y_\eps=Star(\eta(M_{(0,\eps]}))$, see Proposition
\ref{promono}.

In what follows we metrize $X$ so that every $k$-simplex in $X$ is
isometric to the standard Euclidean $k$-simplex; we then will
refer to $X$ as being {\em piecewise-Euclidean}. Our strategy is
to prove an analogue of the vanishing Theorem \ref{vanishing}
first for $(X, Y)$, and then use Proposition \ref{promono} to
derive the desired conclusion for $(M, M_{(0,\eps]})$.

\begin{prop}
\label{nu} There exists a constant $\nu=\nu_m$ with the following
property. Let $X$ be an $m$-dimensional piecewise-Euclidean
simplicial complex, $\W\to X$ be a flat bundle over $X$ and
$Y\subset X$ be a subcomplex. Let $[\zeta]\in H_q(X, Y; \W)$ be a
relative class of dimension $q\ge 1$ so that $Vol([\zeta],
Y)<\nu$. Then $[\zeta]=0$.
\end{prop}
\proof The idea is to retract $\zeta$ inductively to the
$q$-dimensional skeleton of $X$ (away from $Y$) without increasing
the volume too much. The resulting cycle $\zeta'$ will have
relative volume which is less than the volume of the Euclidean
$q$-simplex, therefore $\zeta'$ will miss a point in every
$q$-simplex in $\ol{X\setminus Y}$. Then we retract $\zeta'$ to
the $q-1$-dimensional skeleton of $X$ away from $Y$, thereby
proving vanishing of $[\zeta]$.

\begin{lem}
Suppose that $\Del\subset \ol{X\setminus Y}$ is a $k$-simplex, $k\ge 1$.
There exists a constant $D=D(k)$ such that for every $i< k$ the following holds.

Let $\tau\in C_i(\Del, \W)$. Then there exists a point $x\in \Del$ which avoids the support of $\tau$ and
a retraction $r: \Del\setminus \{x\} \to \partial \Delta$ so that
\begin{equation}
\label{ff}
Vol(r(\tau))\le D\cdot Vol(\tau).
\end{equation}
\end{lem}
\proof This lemma (called {\em Deformation Lemma}) was proved in \cite{Federer} in the case of the
trivial $\R$-bundle $\W$ over $\Del$.
Since our bundle $\W$ is trivial over $\Delta$, the map $r$ defined in \cite{Federer},
extends to the restriction $\W|\Delta$ by the identity along the fibers. Since the volume of a chain is defined
independently of the coefficients, it follows that the inequality (\ref{ff}) holds for the general flat bundles. \qed

We set
$$
D:= \max \{ D(i), q+1\le i\le m\}.
$$

\begin{lem}\label{ffc}
Let $\tau\in C_i(X; \W)$ be a chain. Then there exists another chain $\tau'\in C_i(X, \W)$ so that:

1. The support of $\tau'$ away from $Y$ is contained in the $i$-skeleton of $X$.

2. $Vol(\tau')\le D^{m-i}  Vol(\tau)$.

3. If $\tau\in Z_i(X, Y; \W)$, then
$$
[\tau]=[\tau']\in H_i(X, Y; \W).
$$
\end{lem}
\proof We apply Lemma \ref{ff} inductively. Start with the $m$-skeleton of $X$. For each
$m$-simplex $\Del$ which is not contained in $Y$, we apply the retraction
$$
r: \Delta\setminus \{x\} \to \partial \Delta
$$
adapted to the $i$-chain $\tau\cap \Delta$ obtained from $\tau$ by
excising $X\setminus \Delta$.  We do nothing for the simplices
which are contained in $Y$. The result is an $i$-chain $\tau_1$
such that, away from $Y$, the support of $\tau_1$ is contained in
the $m-1$-skeleton of $X$. By Lemma \ref{ff},
$$
Vol(\tau_1)\le D\cdot Vol(\tau).
$$
We now repeat the above procedure with respect to the $m-1$-skeleton of $X$ and continue inductively $m-i$ times. \qed

\medskip
Let  $v_q$ denote the volume of the standard Euclidean $q$-simplex.
Set
$$
\nu=\nu(q):= v_q\cdot D^{q-m}.
$$
Let $[\zeta]\in H_q(X, Y; \W)$ be such that $Vol(\zeta, Y)<\nu$.
We claim that $[\zeta]=0$. Indeed, by applying Lemma \ref{ffc}, we
construct a relative cycle $\zeta'\in Z_q(X, Y; \W)$ which is
homologous to $\zeta$ and such that
$$
Vol(\zeta')\le   D^{m-q}  Vol(\zeta) <\nu.
$$
Therefore, the support of $\tau$ away from $Y$ is contained in the $q$-skeleton of $X$.
Moreover, for every $q$-simplex $\Del\subset \ol{X\setminus Y}$,
$$
Vol(\zeta'\cap \Delta)\le Vol(\zeta') < v_q=Vol(\Delta).
$$
Therefore $\zeta'$ misses a point $x_\Del$ in the interior of every $q$-simplex $\Del\subset  \ol{X\setminus Y}$.
For every such $q$-simplex $\Del$ we apply the retraction
$$
\rho_\Del: \Del\setminus \{x_\Del\}\to \partial \Del
$$
to the relative cycle $\zeta'$. The result is a new relative cycle $\zeta''$ which is homologous to $\zeta'$
and whose support away from $Y$ is contained in the $q-1$-skeleton of $X$. Since
$$
H_q(X^{(q-1)}, Y^{(q-1)}; \W)=0,
$$
it follows that
$$
[\zeta]=[\zeta']=[\zeta'']=0\in H_q(X, Y; \W).
$$
Lastly, set $\nu_m:=\max \{ \nu(q) : 0< q\le m\}$. Proposition
\ref{nu} follows.  \qed

\medskip
We are now ready to prove Theorem \ref{vanishing}.
Choose $\kappa>0$ so that
$$
\eps+2\kappa< \mu_n,
$$
i.g., take $\kappa=\eps$. 
By Proposition \ref{covering}, there exists a covering $\U=\{B_i, i\in I\}$ of the manifold $M$ by $\rho_i$-balls $B_i$,
where
$$
\frac{\eps}{16}= \sup_{i\in I} \rho_i,
$$
$I\subset \N$. We obtain a piecewise-smooth map
$$
\eta: M\to X,
$$
to the nerve of this covering. The restriction of $\eta$ to $M_{[\eps,\infty)}$ is $L$-Lipschitz,
where $L=\sqrt{m}L_\eps$,  see Corollary \ref{Lip}.
Set
$$
\theta:= \frac{\nu}{L^q},
$$
where $\nu=\nu_m$ is given by the Proposition \ref{nu}.

Consider a cycle $\zeta\in Z_q(M, M_{(0,\eps]}; \V)$. Then
$$
\eta_\#(\zeta\cap M_{(0,\eps]})\subset Y.
$$
Therefore
$$
Vol(\eta_\#(\zeta), Y)\le Vol(\eta_\#(\zeta\cap M_{[\eps,\infty}),
Y)\le L^q Vol(\zeta).
$$
Hence
$$
Vol(\zeta, Y)<\theta \Rightarrow Vol(\eta_\#(\zeta), Y)< \nu \Rightarrow [\eta_\#(\zeta)]=0\in H_q(X, Y_\eps; \W),
$$
by Proposition \ref{nu}. Proposition \ref{promono} implies that
$$
[\zeta]=0\in H_q(M, M_{(0,\eps+2\kappa]}; \V)  \Rightarrow [\zeta]=0\in H_q(M, M_{(0,\mu_n]}; \V),
$$
since $\eps+2\kappa<\mu_n$. Since the map (induced by the inclusion of pairs)
$$
H_q(M, M_{(0,\eps]}; \V)  \to H_q(M, M_{(0,\mu_n]}; \V)
$$
is an isomorphism, it follows that
$$
[\zeta]=0\in H_q(M, M_{(0,\eps]}; \V). \qed
$$

\begin{cor}\label{cor:vanishing}
Let $M$ be a complete (connected) hyperbolic $n$-manifold as above.
Let $\V\to M$ be the flat bundle associated with a $R\Ga$-module $V$.
Then every relative homology class
$$
[\zeta]\in H_q(M, M_{(0,\eps]}; \V), q>0,
$$
of zero relative volume, is trivial.
\end{cor}

\begin{cor}
\label{cuspvanish}
Let $\zeta\in Z^{lf}_p(M; \V)$ be such that $Vol(Exc_{\eps}(\zeta))<\theta(n, \eps)$ for some $0<\eps\le \mu_n/4$.
Then $\zeta$ projects to $0$ in $H^{cusp}_p(M; \V)$.
\end{cor}

\section{Proof of Theorem \ref{main}}

Let $\Ga\subset \Isom(\H^n)$ be a Kleinian group and $\Pi\subset
\Ga$ a  collection of parabolic subgroups as in the Introduction.
Without loss of generality, we may assume that $\Ga$ is
torsion-free. Then the quotient $M=\H^n/\Ga$ is a hyperbolic
manifold.
Set $\del:=\del(\Ga)$. 
Consider an arbitrary $R\Ga$-module $V$ and the relative homology group
$$
H_q(\Ga, \Pi; V),
$$
for $q> \del+1$. Note that $q\ge 2$. Set
\begin{equation}\label{lambda}
\la:= \left( \frac{\del +1}{q}\right)^q.
\end{equation}
Since $q> \del+1$, it follows that $\la<1$.

As it was explained in Section \ref{prelim}, we can use the manifold $M=\H^n/\Ga$  in order to compute the relative
homology of $\Ga$. Let $0<\eps< \mu_n/4$, where $\mu_n$ is the Margulis constant for $\H^n$.
Consider the $\eps$-thick-thin decomposition of the manifold $M$:
$$
M= M_{(0,\eps]}\cup M_{[\eps,\infty)}
$$
and let $K:=M_{(0,\eps]}$. Let $P$ denote the union of components of $K$ whose fundamental group is virtually
abelian of rank $\ge 2$. Then (see Section \ref{prelim})
$$
H_*(M, P; \V)\cong H_*(\Ga, \Pi; V).
$$

\begin{lem}\label{cong}
For every $q\ge 2$,
$$
H_q(M, P; \V)\cong H_q(M, K; \V).
$$
\end{lem}
\proof Set $Q:=K\setminus P$. Then every component of $Q$ has cyclic fundamental group.
Since $hd_R(\Z)=1$ for every ring $R$, it follows that
$$
H_i(Q; \V)=0, i\ge 2.
$$
Therefore for every $q\ge 2$ we have the exact sequence
$$
0=H_q(Q; \V)\to H_q(M, P; \V) \to H_q(M, P\cup Q; \V) \to H_{q+1}(Q; \V)=0.
$$
Hence $H_q(M, P; \V)\cong H_q(M, K; \V)$. \qed

\medskip
Consider a chain $\zeta\in C_q(M; \V)$ which projects to the relative homology class
$$
[\zeta]\in H_q(M, K; \V).
$$
Our goal is to show vanishing of the relative volume of $[\zeta]$:
$$
Vol([\zeta], K)=0,
$$
i.e. that the projection of $\zeta$ to $Z_q(M, K; \V)$ is
homologous to relative cycles of arbitrarily small volume. Then
Theorem \ref{vanishing} will imply that $[\zeta]=0$, thereby
establishing that
$$
0=H_q(M, P; \V)\cong H_q(M, K; \V) \cong H_q(\Ga, \Pi; V),
$$
for all $q>\del+1$.

Let $Q'\subset Q$ be the union of compact components and $P':= K\setminus Q'$.
We start by extending the chain $\zeta\in C_q(M; \V)$ to the tubes in $K$, to a chain
$\zeta'$, which projects to a relative homology class
$$
[\zeta']=[\zeta]\in H_q(M, P'; \V),
$$
see Section \ref{tubes}. We then extend $\zeta'$ to the cusps, to a locally-finite cycle
$$
\hat\zeta\in Z_q(M; \V),
$$
so that
$$
v=Vol(\hat\zeta)\le Vol(\zeta')+ (q-1) Vol(\partial \zeta').
$$
Moreover,
$$
[\Exc_\kappa(\hat\zeta)]=[\zeta']=[\zeta]\in H_q(M, P'; \V)
$$
for every $0<\kappa\le \eps$. See Section \ref{cusps}. We will
need

\begin{thm}
(Besson, Courtois, Gallot, \cite{BCG03, BCG05}) There exists a
smooth map\newline  $\tilde{F}: \H^n\to \H^n$ so that:

\begin{equation}\label{1}
\tilde{F}\circ \ga= \ga\circ \tilde{F}, \quad  \forall \ga\in \Ga.
\end{equation}

\begin{equation}\label{2}
|Jac_r(\t{F}(x))|\le  \left( \frac{\del +1}{r}\right)^r, \quad \forall x\in \H^n, \forall r\ge 1.
\end{equation}

\begin{equation}\label{3}
\t{F}(Hull(\La(\Ga))\subset Hull(\La(\Ga)).
\end{equation}
\end{thm}

The map $\t{F}$ in this theorem is called a {\em natural map}. The
$r$-Jacobian $|Jac_r(\t{F}(x))|$ at $x\in \H^n$ is defined as
$$
\max \{ Vol \left( D_x(\tilde{F})(\xi_1),..., D_x(\tilde{F})(\xi_r) \right) \}
$$
where the maximum is taken over all orthonormal $r$-frames $(\xi_1,...,\xi_r)$ in $T_x\H^n$.
Therefore the map $\tilde{F}$ projects
to a smooth map $F: M\to M$ whose $r$-Jacobian again satisfies (\ref{2})
and which is homotopic to the identity. Hence
$$
Vol(F_\#(\zeta))\le \la Vol(\zeta),
$$
see equation (\ref{lambda}). Since $\la<1$, it follows that
$$
\lim_{k\to\infty} Vol(F^k_\#(\zeta)) =0.
$$
Observe that the inequality (\ref{2}) applied to $r=1$, implies that the map $F$ is $(\del+1)$-Lipschitz.

\smallskip
Let $\theta=\theta_n(\eps)>0$ be as in Theorem \ref{vanishing}.
Since $\la<1$, there exists $k\in \N$ such that
$$
\la^k v<\theta.
$$

Set $f:=F^k$. Then for every $\kappa\le \eps$,
$$
Vol(f_{\#}( \Exc_{\kappa}(\hat\zeta))) \le Vol(f_{\#}(\hat\zeta))\le \la^k v < \theta.
$$

Choose $\kappa:= \frac{\eps}{(1+\del)^k}$ and set
$$
\zeta'':=f_{\#}( \Exc_{\kappa}(\hat\zeta)).
$$

\begin{prop}\label{homologous}
$$
[\zeta'']= [\zeta]
$$
in $H_q(M, M_{(0,\eps]}; \V)$.
\end{prop}
\proof 1. First, we have to check that
$$
\zeta''\in Z_q(M, M_{(0,\eps]}; \V).
$$
Since $f$ is $(\del+1)^k$--Lipschitz, it sends the
$\kappa$--thin part of $M$ to the $(\del+1)^k\kappa$--thin part of $M$.
Therefore
$$
\partial \zeta''\in C_q(M_{(0,\eps]}; \V),
$$
which implies our assertion.

2. Since $f$ is homotopic to the identity and
the cusps $P_i$ are pairwise disjoint, we see that
$$
f(P_i \cap M_{(0,\kappa]} ) \subset P_i
$$
for every component $P_i\subset P'$.

We define the straight-line homotopy $h_t: f\cong Id$ by projecting the straight-line homotopy
$$
\t{f}:=\t{F}^k \cong Id
$$
in $\H^n$. The equality
$$
[\zeta'']= [\zeta]
$$
would follow from

\begin{lem}
For every $x\in P_i \cap M_{(0,\kappa]}$, the geodesic
$$
h_t(x)=\ol{x f(x)}
$$
is contained in $P_i$.
\end{lem}
\proof Let $\Pi_i\subset \Ga$  be the fundamental group of $P_i$, i.e.
$\Pi_i$ is the stabilizer in $\Ga$ of a component $\t{P}_i$ of the lift of $P_i$ to $\H^n$.
Recall that $P_i \cap M_{(0,\kappa]}$ is the projection to $M$ of the union
$$
\bigcup_{\ga\in \Pi_i\setminus \{1\} } K_\kappa(\ga),
$$
see Section \ref{prelim}. Since $\t{f}=\t{F}^k$ is $(1+\del)^k$-Lipschitz and commutes with $\Ga$, we obtain
$$
f(K_\kappa(\ga))\subset K_{(1+\del)^k\kappa} (\ga)\subset K_{\eps}
(\ga),
$$
for all $\ga\in \Pi_i\setminus \{1\}$. Since $K_{\eps} (\ga)$ is convex, for every $\t{x}\in K_\kappa(\ga)$,
$$
\ol{ \t{x} \t{f}(z) } \subset  K_{\eps} (\ga).
$$
The above geodesic segment projects to the track $\ol{x f(x)}$ of the homotopy $h_t$ connecting $x=p(\t{x})$ to $f(x)$.
On the other hand, $K_{\eps} (\ga)$ projects to $P_i$. Therefore
$$
\ol{x f(x)}\subset P_i. \qed
$$

This concludes the proof of Proposition \ref{homologous}. \qed

\medskip
We now can finish the proof of Theorem \ref{main}. Since
$$
Vol(\zeta'')<\theta=\theta_n(\eps)
$$
and $0<\eps\le \mu_n/4$, Theorem \ref{vanishing} implies that $[\zeta'']=0$ in
$$
H_q(M, M_{(0,\eps]}; \V).
$$
By combining this with Proposition \ref{homologous}, we obtain
$$
[\zeta]=[\zeta'']=0\in H_q(M, M_{(0,\eps]}; \V).
$$

Therefore, for every $p>\del+1$,
$$
0= H_q(M, P; \V)\cong H_q(\Ga, \Pi; V).
$$
Hence
$$
vhd_R(\Ga)\le \del(\Ga)+1.
$$
This concludes the proof of Theorem \ref{main}. \qed

\bigskip
{\em Proof of} Corollary \ref{maincor}. Since $(\Ga, \Pi)$ has
finite type, it follows that
$$
cd_R(\Ga, \Pi)=hd_R(\Ga, \Pi)=vhd_R(\Ga, \Pi)\le \del(\Ga)+1. \qed
$$

\bigskip
{\em Proof of} Corollary \ref{free}. Since $\Ga$ is of type
$FP_2$, it is finitely-generated, hence $\Ga$ is virtually
torsion-free by Selberg's lemma. Therefore, without loss of
generality we may assume that $\Ga$ is torsion-free. Since
$\del(\Z^2)=1$, it follows that $\Ga$ contains no free abelian
subgroups of rank $\ge 2$. Thus $\Pi=\emptyset$. By Theorem
\ref{main}, we have the inequalities
$$
cd_R(\Ga)\le 1+ hd_R(\Ga)\le \delta(\Ga)+2<3.
$$
The above inequality implies that $cd_R(\Ga)\le 2$. Applying Lemma
\ref{brown} to the group $\Ga$, we conclude that $\Ga$ is of type
$FP$ and hence
$$
cd_R(\Ga)=hd_R(\Ga),
$$
see Lemma \ref{L2}.  Applying Theorem \ref{main} again, we obtain
the inequality
$$
cd_R(\Ga)=hd_R(\Ga)\le \del(\Ga)+1<2.
$$
Therefore $cd_R(\Ga)=1$ and hence $\Ga$ is free by Theorem
\ref{Stallings-Dunwoody}. \qed

\begin{rem}
If one could replace $vhd$ with $vcd$ in Theorem \ref{main}, then
one can weaken the assumption in Corollary \ref{free} to finite
generation of $\Ga$.
\end{rem}

\section{Application to geometrically finite groups}

The main goal of this section is to prove Theorem \ref{lattice}
from the introduction. In what follows, let $\Ga\subset
\Isom(\H^n)$ be a nonelementary geometrically finite Kleinian
group. Without loss of generality we may assume that $\Ga$ is
torsion-free. Recall that $\Pi$ is a maximal collection of maximal
parabolic subgroups of $\Ga$ (of virtual rank $\ge 2$) which are
pairwise nonconjugate  and such that every maximal parabolic
subgroup of $\Ga$ is conjugate to one of the subgroups $\Pi_i$.

We enlarge $\Pi$ to the set
$$
\Pi'=\{\Pi_i, i\in I\},
$$
which consists of representatives of conjugacy classes of {\em
all} maximal parabolic subgroups in $\Ga$. Note that
$$
hd_R(\Ga, \Pi)=hd_R(\Ga, \Pi'),
$$
see the proof of Lemma \ref{cong}.

We will need the following

\begin{prop}
Let $\Ga\subset \Isom(\H^n)$ be a discrete subgroup and $\t{F}: \H^n\to \H^n$ be the natural map associated with $\Ga$.
Suppose that there exists $\t{x}\in \H^n$ such that $|Jac_q(\t{F}(\t{x}))|=1$. Then there exists
a $q$-dimensional subspace $H\subset \H^n$ through $x$, so that $Hull(\La(\Ga))$ is contained in $H$.
\end{prop}
\proof This proposition was proved in \cite[Proposition 5.1]{BCG05} in the case $q=n-1$. It is clear
from their proof however that it works for arbitrary $q$. \qed

We will assume in what follows that $\Ga$ does not preserve any
proper subspace in $\H^n$ (otherwise we pass to the smallest
$\Ga$-invariant subspace). Then the convex hull $Hull(\La(\Ga))$
of $\La(\Ga)$ is $n$-dimensional.

\begin{cor}
\label{latco} If there exists $x\in M$ such that
$|Jac_q(F(x))|=1$, then $q=n$.
\end{cor}

\noindent The key technical result of this section is the
following strengthening of Theorem \ref{main}:

\begin{prop}
\label{subspace} Suppose that $\Ga\subset \Isom(\H^n)$ is a
geometrically finite group such that $\del(\Ga)+1=cd_R(\Ga,
\Pi')=q$. Then $q=n$ and $\La(\Ga)=S^{n-1}$.
\end{prop}
\proof Set $\del:=\del(\Ga)$. Let $N=Hull(\La(\Ga))/\Ga$ denote
the convex core of $M$; then $N$ is $n$-dimensional. Let $\t{F}:
\H^n\to \H^n$ be the natural map associated with $\Ga$. Our goal
is to find a point $\t{x}\in \H^n$ such that
$|Jac_q(\t{F}(\t{x}))|=1$. Once we found such a point, it will
follow from Corollary \ref{latco} that $q=n$ and hence
$$
\del=q-1=n-1.
$$
Then, since $\Ga$ is geometrically finite, it will follow that
$\La(\Ga)=S^{n-1}$, see Theorem \ref{nich}, Part 3.

Note that the projection $F: M\to M$ of the map $\t{F}$ satisfies
$$
F(N)\subset N.
$$
Moreover, since $\del(\Ga)+1=q$, the inequality (\ref{I2}) implies
that $F$ does not increase the volume of chains in $C_q^{lf}(M;
\V)$. Since $F$ is $(1+\del)$-Lipschitz, the restriction $F|N$ is
a proper map.

Since $N=K(\Ga,1)$, we can use the thick part of this manifold in
order to compute the (relative) homology of $\Ga$: We choose
$0<\eps<\mu_n/4$ such that
$$
N_{(0,\eps]}=P_\eps,
$$
is the disjoint union of cusps. Since $\Ga$ is geometrically
finite, it follows that $N_{[\eps,\infty)}$ is compact. Since
$q=cd(\Ga, \Pi')$ and the pair $(\Ga,\Pi')$ has finite type,
$$
q=cd_R(\Ga,\Pi')=hd_R(\Ga,\Pi').
$$
Hence there exists  an $R\Ga$-module $V$ and a nonzero relative
homology class
$$
[\zeta]\in H_q(N, P_\eps; \V)\cong H_q(\Ga, \Pi'; V).
$$
Let
$$
\hat\zeta:=\ext(\zeta)\in Z^{lf}_q(N; \V)$$
be a finite volume extension of the relative cycle $\zeta$.
Since $Vol(\hat\zeta)<\infty$, there exists $0<\kappa<\eps$ such that
the chain $\zeta'':=\hat\zeta\cap P_\kappa$ satisfies
$$
Vol(\zeta'')<t:=\theta_n(\eps)/2,
$$
where $\theta_n(\eps)$ is the function introduced in Theorem
\ref{vanishing}. Therefore
$$
Vol(F^k(\zeta''))<t,
$$
for all $k\ge 0$. Set $\zeta':= Exc_\kappa(\hat\zeta)=
\hat\zeta\cap N_{[\kappa,\infty)}$.

For $k\in \N$ define the chains
$$
\zeta'_k:= F^k(\zeta'), \quad \zeta_k^+:=\zeta_k'\cap
N_{[\eps,\infty)}, \quad \zeta_k^-:= \zeta_k'\cap N_{(0, \eps]}.
$$

We will consider two cases, in the first case we find a point
$y\in N$ such that $|Jac_q(F(y))|=1$, the second case will be
ruled out as it will lead to the contradiction with nonvanishing
of $[\zeta]$.

\smallskip
{\bf Case 1.} Let $\chi$ denote the characteristic function of $N_{[\eps,\infty)}$. Suppose that
there exists a sequence $x_j\in Supp(\zeta')$ and $k_j\in \N$, such that for $y_j:=F^{k_j}(x_j)$ we have:
$$
\lim_j  \chi(y_j)\cdot |Jac_q(F( y_j))| =1.
$$
Then the sequence $(y_j)$ belongs to the compact $N_{[\eps,\infty)}$ and hence subconverges to
a point $y$ so that
$$
|Jac_q(F(y))|=1.
$$
Thus we are done by Corollary \ref{latco}.

\smallskip
{\bf Case 2.} Otherwise, there exists $0<\la<1$ so that for all $k\in \N$,
\begin{equation}
\label{la}
Supp(\zeta_k^+)\subset E_\la:=\{x\in M: |Jac_q(F(x))|\le \la\}.
\end{equation}

\begin{lem}
There exists $k\in \N$ such that $Vol(\zeta'_k, P_\eps)<
t=\theta_n(\eps)/2$.
\end{lem}
\proof Suppose not. Then for every $k$ we have
$$
Vol(\zeta^+_k) \ge t
$$
and hence
$$
\frac{Vol(\zeta^+_k)}{Vol(\zeta'_k)}\ge \frac{t}{v},
$$
where
$$
v:= Vol(\zeta')\ge Vol(\zeta'_k).
$$
Moreover, by (\ref{la}), we get
$$
Vol(\zeta_{k+1}')=Vol(F(\zeta'_k))\le \la Vol(\zeta_k^+) + Vol(\zeta_k^-)\le [(\la-1)\frac{t}{v}+1]Vol(\zeta'_k).
$$
Note that, since $\la<1$, we have
$$
0< [(\la-1)\frac{t}{v}+1]<1.
$$
Therefore
$$
\lim_{k\to\infty} Vol(\zeta'_k, P_\eps)\le \lim_{k\to\infty} Vol(\zeta'_k)=0.
$$
This contradicts the assumption that $Vol(\zeta'_k, P_\eps)\ge
t>0$ for all $k$. Contradiction. \qed

\medskip
We now can finish the proof of the proposition.
We first estimate $Vol(F^k(\hat\zeta), P_\eps)$ for the number $k$ guaranteed by the
above lemma:
$$
Vol(F^k(\hat\zeta), P_\eps)\le Vol(F^k(\zeta'), P_\eps)+
Vol(F^k(\zeta''))\le t+t=\theta_n(\eps),
$$
since $Vol(F^k(\zeta''))\le t$.

Hence, by Corollary \ref{cuspvanish}, the locally finite cycle $F^k(\hat\zeta)$
projects to zero class
$$
[Exc_\eps(F^k(\hat\zeta))]\in H^{cusp}_q(M; \V).
$$
Since $F^k: N\to N$ is Lipschitz and commutes with $\Ga$, it is
properly homotopic to the identity. Therefore
$$
[\hat\zeta]= F^k([\hat\zeta])\in H^{lf}_q(M; \V).
$$
Thus
$$
[\zeta]=[Exc_\eps(\hat\zeta)]=0\in H_q(M, P_\eps; \V).$$
This contradicts the assumption that $[\zeta]$ is a nonzero class in $H_q(M, P_\eps; \V)$.
\qed

In order to relate the above proposition to the limit set of $\Ga$
we will need the following proposition, which is a relative
version of a theorem by Bestvina and Mess in
\cite{Bestvina-Mess(1991)}:

\begin{prop}\label{BM}
$cd(\Ga, \Pi')=\dim(\La(\Ga))+1$, where $\dim$ is the topological
dimension.
\end{prop}
\proof Without loss of generality we may assume that $\Ga$ is torsion-free.
Let $H:=Hull(\La)$ denote the convex hull of the limit set $\La$ of the group $\Ga$.
The set $H$ is obviously contractible. Moreover, the union
$$
H\cup \La
$$
satisfies the axioms of the ${\mathcal Z}$-set compactification of
$H$ and therefore
$$
H_c^*(H) \cong \check{H}^{*-1}(\La),
$$
see \cite{Bestvina-Mess(1991)} and \cite{Bestvina(1996)}. Since
$\Ga$ is geometrically finite, the pair $(\Ga, \Pi')$ has finite
type;  hence
$$
cd(\Ga, \Pi')= \sup \{ q: H^q(\Ga, \Pi'; \Z\Ga)\ne 0\},
$$
see Lemma \ref{L2}.

Since $\Ga$ is geometrically finite, for every $i\in I$ we can
choose a closed horoball $B_i$, centered at the fixed point of
$\Pi_i$, so that:
$$
\ga(B_i)\cap B_j=\emptyset
$$
unless $i=j$ and $\ga\in \Pi_i$, in which case $\ga(B_i)=B_j$. See
for instance \cite{Bowditch(1993b)}.

Therefore the set
$$
H':=H\setminus \Ga\cdot \bigcup_{i\in I} int(B_i)
$$
projects to a compact submanifold with boundary $N'$ in $N=H/\Ga$.

For every $i\in I$, set
$$
C_i:=B_i\cap H,
$$
$$
C:=\Ga\cdot \bigcup_{i\in I} C_i,
$$
and
$$
C_i':=B_i\cap H',
$$
$$
C':=\Ga\cdot \bigcup_{i\in I} C_i'.
$$

Then convexity of $H$ and of every $B_i$ implies that $H'$ and
each $C'_i$ is contractible. Therefore $N'=H'/\Ga$ is a compact
$K(\Ga, 1)$ and $C_i'/\Pi_i$ is a compact $K(\Pi_i, 1)$, for every
$i\in I$.

\begin{lem}
$H^*_c(C)=0$.
\end{lem}
\proof Since each $B_i$ is a horoball centered at a limit point of
$\Ga$ and $H$ is convex, it follows that
$$
C\cong [0,1)\times C'.
$$
Therefore vanishing of $H^*_c([0,1))$ implies vanishing of $H^*_c(C)$. \qed

Hence, by the long exact sequence of the pair $(H, C)$, we have
$$
H^*_c(H)\cong H^*_c(H, C).
$$
We claim that
$$
H^*(\Ga, \Pi'; \Z\Ga)\cong H^*_c(H', C') \cong H^*_c(H, C)\cong
H^*_c(H)\cong \check{H}^{*-1}(\La).
$$
The first isomorphism in this sequence is established in
\cite[Lemma 7.4]{Brown(1982)} in the case $\Pi=\emptyset$ and
$\Ga$ of finite type with $H'$ being the universal cover of a
compact $K(\Ga,1)$; the general case follows from the long exact
sequences of the pairs $(\Ga, \Pi)$, $(H', C')$. The rest of the
isomorphisms were established above.
 Hence $cd(\Ga, \Pi')=\dim(\La(\Ga))+1$ and Proposition \ref{BM} follows. \qed

\begin{rem}
The proof of Proposition \ref{BM} generalizes without much
difficulty to the case of  $\Ga$ relatively hyperbolic groups with
respect to a family $\Pi'$ of virtually nilpotent subgroups. The
limit set in this case is replaced by the Bowditch boundary of
$\Ga$.
\end{rem}

\medskip
We are now ready to prove Theorem \ref{lattice}. By the assumption,
$$
\dim_H\La(\Ga)=\dim \La(\Ga).
$$
Since, by Theorem \ref{nich}, for geometrically finite Kleinian
groups $\Ga$ we have
$$
\dim_H\La(\Ga)=\del(\Ga),
$$
Proposition \ref{BM} implies
$$
\del(\Ga)= \dim_H\La(\Ga)=\dim \La(\Ga)=cd(\Ga, \Pi').
$$
Lastly, Proposition \ref{subspace} implies that $\Ga$ is a lattice. \qed

\bibliography{lit}
\bibliographystyle{siam}

\medskip

Department of Mathematics,

University of California,

Davis, CA 95616, USA,

kapovich$@$math.ucdavis.edu

\end{document}